\documentclass[conference]{IEEEtran}
\IEEEoverridecommandlockouts
\usepackage{cite}
\usepackage{amsthm}
\usepackage{amsmath,amssymb,amsfonts}
\usepackage{algorithmic}
\usepackage{import}
\usepackage{graphicx}
\usepackage{textcomp}
\usepackage{xcolor}
\usepackage{subfigure}
\usepackage{url} 
\usepackage{geometry}

\usepackage{xcolor,calc}

\newgeometry{left=1.91cm, right=1.91cm, top=2.54cm, bottom=1.91cm}

\newtheoremstyle{spaceless} 
  {0pt} 
  {0pt} 
  {\itshape} 
  {} 
  {\bfseries} 
  {.} 
  { } 
  {} 

\theoremstyle{spaceless}
\newtheorem{lemma}{Lemma}[section]
\newtheorem{definition}{Definition}[section]
\newtheorem{theorem}{Theorem}

\def\BibTeX{{\rm B\kern-.05em{\sc i\kern-.025em b}\kern-.08em
    T\kern-.1667em\lower.7ex\hbox{E}\kern-.125emX}}
\begin{document}

\title{\LARGE \bf
Predictive Control Barrier Functions: Bridging model predictive control and control barrier functions\\
}

\author{Jingyi Huang, Han Wang, Kostas Margellos, Paul Goulart

\thanks{The authors are with the Department of Engineering Science,
        University of Oxford, Oxford, United Kingdom. E-mails:
        \{jingyi.huang, han.wang, kostas.margellos, paul.goulart\}@eng.ox.ac.uk}%

\thanks{For the purpose of Open Access, the authors have applied a CC BY copyright licence to any Author Accepted Manuscript (AAM) version arising from this submission.}
}

\maketitle

\begin{abstract}
In this paper, we establish a connection between model predictive control (MPC) techniques and Control Barrier Functions (CBFs). Recognizing the similarity between CBFs and Control Lyapunov Functions (CLFs), we propose a MPC formulation that ensures invariance and safety without relying on explicit stability conditions. The value function of our proposed MPC is a CBF, which we refer to as the Predictive Control Barrier Function (PCBF), similar to traditional MPC formulations which encode stability by having value functions as CLFs. Our formulation is simpler than previous PCBF approaches and is based on weaker assumptions while proving a similar theorem that guarantees safety recovery. Notably, our MPC formulation does not require the value function to be strictly decreasing to ensure convergence to a safe invariant set. Numerical examples demonstrate the effectiveness of our approach in guaranteeing safety and constructing non-conservative CBFs.
\end{abstract}

\section{Introduction}
The development of new control methods that can provide safety guarantees is critical due to the increasing prevalence of autonomy in safety-critical applications. The challenge of safety is that systems under unsafe control actions may not immediately violate state constraints, but can nevertheless inevitably reach unsafe regions. Safety can therefore be guaranteed by ensuring positive invariance of the safe set \cite{prajna2005necessity}. A set is  positive invariant if the inclusion of the state at some times implies the inclusion in the future \cite{blanchini1999invariant_set}.

One common way to verify the positive invariance of a safe set is through reachability analysis, which computes the set of states that can reach a target set while remaining in the safe set \cite{ Recent_reachability, HJ_formulation, HJ_kostas}. Reachability-based methods construct the \emph{maximal safe invariant set}, but they are well-known to suffer from the curse of dimensionality \cite{Reachability_overview}. 
To address this, numerical methods such as SOS \cite{SOS_formulation, Han_book}, learning \cite{robey2020learningcontrolbarrierfunctions, Learning_certificate_function}, etc., have been proposed to approximate the value function of the reachability optimal control problem. The approximate function, referred to as \emph{Control Barrier Function (CBF)} \cite{Alan2022ControlBF, CBF_theory_and_application, CDC_NPM_2023}, identifies a (potentially conservative) invariant subset of the safe set. Current methods usually yield overly conservative CBFs, and efficient construction of guaranteed and non-conservative CBFs remains an open challenge.

The conditions that determine whether a given function is a CBF are strongly reminiscent of those defining \emph{Control Lyapunov Functions (CLFs)}. However, unlike CLFs, which require invariance for \emph{every} sublevel set, CBFs require invariance only of their zero sublevel set. 
Previous research on MPC \cite{MPC_optimality, MPC_book} has shown that MPC can ensure stability by proving its value function is a CLF, extending a local control Lyapunov function defined over the terminal set. Given the similarities between the governing equations of CBFs and CLFs, we establish a similar connection between MPC and CBFs, in particular showing that MPC formulation can guarantee safety, with its value function acting as a CBF, referred to as a \emph{Predictive Control Barrier Function (PCBF)}, that increasingly enlarges a local conservative invariant set with longer horizons. 

In this work we propose a safe MPC formulation to guarantee positive invariance and safety. It constructs a non-conservative safe invariant set and ensures convergence of trajectories that start from initially unsafe regions to the maximal safe invariant set. In contrast to traditional MPC convergence proofs that rely on Lyapunov's theorem \cite{MPC_optimality}, requiring strictly decreasing value functions, we apply the local invariant set theorem \cite[Section 3.4.3]{Nonlinear_Control_Book} to guarantee convergence with weaker conditions on a non-increasing value function. 

The authors of \cite{Predictive_CBF} proposed a similar soft-constrained MPC formulation, but it further tightens the soft state constraints to have a strictly decreasing value function to guarantee convergence.  We show that the tightening is not necessary to ensure convergence of the safe invariant set, and therefore to recover safety.

The idea of using a receding horizon to enlarge an invariant set has also been explored in \cite{backup_CBF}, where a nominal trajectory is propagated using a predetermined backup control policy $\pi$. However, the size of the enlarged invariant set depends heavily on the chosen backup policy, while our formulation does not require such a policy.

\textit{Contributions:}
The main contributions of this paper are:
\begin{itemize}
    \item We propose a safe MPC formulation that guarantees positive invariance and convergence towards the maximal safe invariant set.
    \item We prove that the value function of the MPC problem acts as a control barrier function (CBF), in a manner analogous to the development of control Lyapunov functions in standard MPC methods.  
    \item We apply our MPC formulation to linear and nonlinear examples, demonstrating its ability to expand a predef-\restoregeometry ined conservative safe invariant set and guarantee local convergence of trajectories to a safe invariant set.
\end{itemize}

\textit{Outline:} The paper is organized as follows: Section II formulates the problem and presents the necessary preliminaries. Section III proposes the safe MPC formulation and proves its convergence and stability properties. It also shows that the value function of the safe MPC problem is a CBF, providing a non-conservative inner approximation of the maximal safe invariant set. Section IV shows two numerical examples to illustrate a safety recovery mechanism and compare our method with the literature. Section V concludes the paper.

\textit{Notation:} The distance from a point $x \in \mathbb{R}^n$ to a set $A \subseteq \mathbb{R}^n$ is defined as $|x|_A := \inf_{y \in A}||x-y||$, and $||\cdot||$ denotes the Euclidean norm. Given a closed set $C$, $\partial C$ and $\text{Int}(C)$ are the boundary and the interior of the set $C$. The notation $D \backslash S$ represents the relative complement of $S$ in $D$. In the MPC formulation, closed-loop elements are denoted with brackets, e.g., $x(k)$ is the closed-loop state at time $k$. A symbol with two subscripts represents an open-loop element, where the first subscript indicates the lookahead time and the second subscript indicates the current time, e.g., $x_{i|k}$ denotes the open-loop state predicted $i$ steps ahead when evaluating at time $k$. The optimal solution has superscript $*$, e.g., $\mathbf{x}_k^*$ is the optimal open-loop state sequence at current time $k$. A proposed solution is denoted with a hat, e.g., $\hat{\mathbf{x}}_k$ is the proposed open-loop state sequence at current time $k$. This paper considers a time-invariant system, so the MPC optimal solution depends only on the current system state $x$, to make this dependency explicit, we denote the closed-loop optimal slack variable and closed-loop control input at state $x$ as $\xi_{0}^*(x)$ and $u_0^*(x)$, respectively.

\section{Preliminaries}
In this section we formulate the problem of interest and present the necessary preliminaries.
\subsection{Problem formulation}
We consider a discrete-time system
\begin{equation}\label{system_dynamics}
    x(k+1) = f(x(k), u(k)), \; k \in \mathbb{N},
\end{equation}
with continuous dynamics $f: \mathbb{R}^n \times \mathbb{R}^m \rightarrow \mathbb{R}^n$ and initial condition $x(0) \in \mathbb{R}^n$. The system is subject to hard physical input constraints of the form 
\begin{equation}\label{input_constraint
}
    u(k) \in \mathbb{U} \subset \mathbb{R}^m, \; k \in \mathbb{N}, 
\end{equation}
where $\mathbb{U}$ is a polytopic set. The system is required to satisfy safety constraints that are formulated as state constraints given by 
\begin{equation}\label{state_constraint}
    x(k) \in \mathbb{X} \subset \mathbb{R}^n, \; k \in \mathbb{N},
\end{equation}
along the trajectory. We assume that $\mathbb{X}$ is a polytopic set, represented as $\mathbb{X} := \{x \in \mathbb{R}^n : c(x) \leq 0\}$, where $c: \mathbb{R}^n \rightarrow \mathbb{R}^{n_x}$. In this work, we focus on guaranteeing safety and set stability rather than stability around an equilibrium point. Thus, our goal is to construct a control law $u(k) = \kappa(x(k))$ that finds the maximal invariant subset of $\mathbb{X}$ and steers the system towards it even from some initially unsafe states, the latter can be thought of as a safety recovery mechanism. 

\subsection{Safety and barrier certificates}
Barrier functions constitute a class of certificate functions, used to prove positive invariance of a set and ensure safety.

Consider a safe set $C$, which can be defined as 
\begin{equation}\label{C}
\begin{split}
    C & := \{x\in \mathbb{R}^n: h(x) \leq 0 \}, \\
    \partial C& := \{x \in \mathbb{R}^n: h(x) = 0 \}, 
\end{split}
\end{equation}
where $h$ is a continuous function. The set $C$ is assumed to have a non-empty interior, i.e., $\text{Int}(C) \neq \emptyset$. 

\begin{definition}[Discrete-time Control Barrier Function]
Consider the safe set $C$ defined by \eqref{C}. If there exist a continuous function $h:\mathbb{R}^n \rightarrow \mathbb{R}$, defined on a set $D$ with $C \subseteq D \subset \mathbb{R}^n$, such that 
\begin{equation}\label{self_defined_h_condition}
\inf_{u\in\mathbf{U}} [h(f(x,u)) - h(x)] \leq 0, \forall x \in D,
\end{equation}
then the function $h$ is a \emph{Discrete-Time Control Barrier Function} with domain $D$.
\end{definition}

If such a Discrete-Time CBF \( h \) exists, then the safe set \( C \) is positively invariant, ensuring safety. Specifically, for all $x \in C,\; h(x) \leq 0$, according to \eqref{self_defined_h_condition} there exists a controller \( \kappa_{CBF}(x) = \{ u \in \mathbf{U} : h(f(x,u)) - h(x) \leq 0 \} \). This implies \( h(f(x,u)) \leq h(x) \leq 0 \) and \( f(x,u) \in C \), proving that \( C \) is positively invariant.


\subsection{Stability and convergence results}
In order to establish our results in Section III, we will require the following preliminaries:
\begin{lemma}[\protect{{\cite[Lemma 3.4]{bertsekas1996neuro}}}]
\label{convergence}
Consider non-negative scalar sequences $(l(k))_{k \in \mathbf{N}}, (w(k))_{k \in \mathbf{N}}$ and $(\zeta(k))_{k \in \mathbf{N}}$  that satisfy the recursion $l(k+1) \leq l(k) - w(k) + \zeta(k)$. If $\sum^{\infty}_{k=0} \zeta(k) < \infty$, then the sequence $(l(k))_{k\in \mathbf{N}}$ converges and $\sum^{\infty}_{k=0} w(k) < \infty$.
\end{lemma}

Standard stability theory considers the stability of an equilibrium point and can be extended to the stability of a set. In \cite{Predictive_CBF}, the concept of asymptotic stability for an invariant set is established, and we adapt this definition to the stability of an invariant set.

Consider system \eqref{system_dynamics} using the closed-loop control law $u(k) = \kappa(x(k))$, i.e., 
\begin{equation}\label{closed_loop_system}
    x(k+1) = f(x(k), \kappa(x(k))), \; k \in \mathbb{N},
\end{equation}
subject to input and state constraints, $\mathbb{U}$ and $\mathbb{X}$, with initial state $x(0)$. The stability of an invariant set with respect to the system \eqref{closed_loop_system} can be defined as follows.

\begin{definition}[Stability with respect to a set, {\protect{\cite[Definition III.3]{Predictive_CBF}}}] Let $S$ and $D$ be closed, non-empty and positively invariant sets for the closed-loop system \eqref{closed_loop_system}, with $S \subset D$. The set $S$ is a stable set in $D$, if for all $x(0) \in D$ the following conditions hold:
\begin{equation}\label{stability_definition}
    \forall \epsilon > 0, \exists \delta >0: |x(0)|_S<\delta \Rightarrow |x(k)|_S < \epsilon, \forall k>0.
\end{equation}
\end{definition}

\begin{definition}[Locally positive definite (l.p.d.) function with respect to sets, {\protect{\cite[Definition A.1]{Predictive_CBF}}}]\label{lpd}
Given two non-empty and closed sets $S$ and $D$, with $S \subset D$, consider a continuous function $V:D\rightarrow \mathbb{R}$. $V(x)$ is called a locally positive definite (l.p.d.) function around $S$ in $D$ if holds that 
\begin{equation*}
\begin{split}
    & V(x) = 0, \forall x \in S \text{ and}\\
    & V(x) > 0, \forall x \in D \backslash S. 
\end{split}
\end{equation*}
\end{definition}

\section{Safe MPC}
\begin{figure}[ht]
\centering
\includegraphics[width=0.35\textwidth]{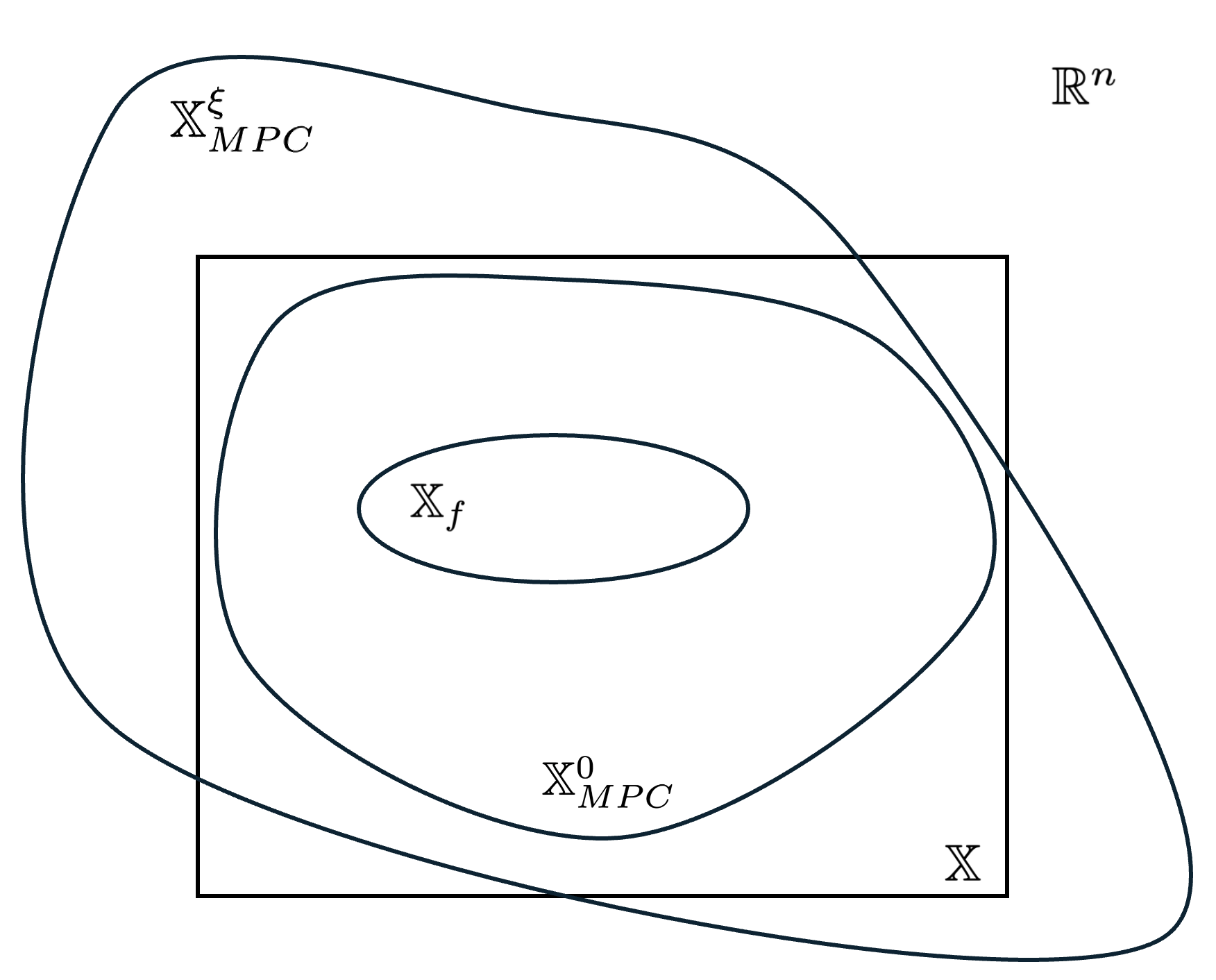}
\caption{The relevant sets in the MPC formulation include: $\mathbb{X}$, the unrelaxed state constraint set; $\mathbb{X}_f$, the terminal set; $\mathbb{X}^{\xi}_{MPC}$, the feasible set for the safe MPC problem \eqref{soft_MPC}; and $\mathbb{X}_{MPC}^0$, the feasible set for the MPC problem without relaxation.}
\vspace{-15pt}
\label{fig:sets}
\end{figure}

We next present our safe MPC formulation with soft state constraints\footnote{While our formulation relaxes all state constraints for simplicity, the theoretical guarantees extend to cases where only partial state constraint relaxation is applied, with the remaining state constraints enforced as hard constraints over the prediction horizon.}. 
We show that all trajectories, starting from the feasible set of the safe MPC formulation, remain in or converge to the maximal safe invariant set of the system. The safe MPC problem is formulated as 
\begin{equation}\label{soft_MPC}
\begin{split}
V^*(x(k)) = \min_{u_{i|k}, x_{i|k}, \xi_{i|k}} & \;\; \sum_{i=0}^{N-1} \xi_{i|k}\\
\text{subject to } & \forall i = 0, 1,...,N-1:\\
& x_{0|k} = x(k), \\
& x_{i+1|k} = f(x_{i|k}, u_{i|k}),\\
& u_{i|k} \in \mathbb{U},\\
& x_{i|k} \in \mathbb{X}(\xi_{i|k}), \xi_{i|k} \geq 0,\\
& x_{N|k} \in \mathbb{X}_f,
\end{split}
\end{equation}
where $\xi_{i|k} \in \mathbb{R}_{\geq 0}$ is a slack variable for the state constraints $i$ steps ahead when evaluated at the current time $k$, $\mathbb{X}(\xi_{i|k})$ is the relaxed state constraint set, where $\mathbb{X}(\xi_{i|k}) := \{x \in \mathbb{R}^n:c(x) \leq \xi_{i|k}\}$\footnote{To relax state constraints elementwise, we may alternatively employ a slack variable vector $\xi_{i|k} \in \mathbb{R}^{n_x}_{\geq 0}$, defining the constrained set as $\mathbb{X}(\xi_{i|k}) := \{x \in \mathbb{R}^n:c(x) \leq \xi_{i|k}\}$. This would lead to an objective function of the form $\sum_{i=0}^{N-1} ||\xi_{i|k}||$. For notational simplicity, we maintain our current formulation, which effectively corresponds to using a slack vector with the $l_\infty$-norm penalty. All theoretical results remain valid when using slack vector formulations.}. 
The terminal set $\mathbb{X}_f$ is compact and $\mathbb{X}_f \subset \mathbb{X}$. We construct $\mathbb{X}_f$ such that there exists a local controller $\kappa_f(\cdot): \mathbb{X}_f \rightarrow \mathbb{U}$, $\forall x \in \mathbb{X}_f$, which ensures the positively invariance of $\mathbb{X}_f$, i.e., $f(x, \kappa_f(x)) \in \mathbb{X}_f, \forall x \in \mathbb{X}_f$. 
We use $\mathbb{X}^{\xi}_{MPC}$ to denote the feasible set of \eqref{soft_MPC}. The feasible set of the same MPC problem without constraint relaxation, i.e., let $\xi_{i|k} = 0$ for all $i$, is denoted by $\mathbb{X}_{MPC}^0$. We assume $\mathbb{X}_{MPC}^0$ to have a non-empty interior. The relationship between the different sets is shown in Fig \ref{fig:sets}.

\textbf{Remarks.} The authors of \cite{Predictive_CBF} proposed a similar formulation, but further tightened every soft state constraint to obtain a strictly decreasing value function and therefore guarantee convergence to the safe invariant set, i.e., recover safety, but in our paper, we show that the tightening is not necessary to converge to the safe invariant set.

\subsection{Convergence of the MPC formulation}
In this section we show that all closed-loop trajectories starting from $ \mathbb{X}^{\xi}_{MPC}$ with the closed-loop safe MPC controller will converge to the safe invariant set within $\mathbb{X} \cap \mathbb{X}^{\xi}_{MPC}$, even for some unsafe initial states $x(0) \notin \mathbb{X}$.

\begin{lemma}\label{lemma_positive_invariant}
Problem \eqref{soft_MPC} defines a controller \( u_0^*(x) \). If \( \kappa(x) = u_0^*(x) \) in system \eqref{closed_loop_system}, then the set \( \mathbb{X}_{MPC}^0 \) and $\mathbb{X}^{\xi}_{MPC}$ are positively invariant.

\end{lemma}

\textit{Proof.} If $x(0) \in \mathbb{X}_{MPC}^0$, i.e., the initial state is feasible for the hard-constrained MPC problem, by the recursive feasibility property \cite{MPC_optimality}, the closed-loop state $x(k)$ of system \eqref{closed_loop_system} is feasible for the same problem, i.e., $x(k) \in \mathbb{X}_{MPC}^0$, for all $k \geq 0$. Thus, $\mathbb{X}_{MPC}^0$ is positively invariant. Similar arguments also hold for the feasible set $\mathbb{X}^{\xi}_{MPC}$ of the safe MPC problem \eqref{soft_MPC}, by noticing that $\mathbb{X}_f$ is positively invariant, which guarantees recursive feasibility of \eqref{soft_MPC}. \qed 

\begin{lemma}\label{V_continuous}
The value function $V^*(\cdot)$ in \eqref{soft_MPC} is continuous over $\mathbb{X}_{MPC}^{\xi}$.
\end{lemma}
\textit{Proof.} Since $\mathbb{X}$ and $\mathbb{U}$ are polytopic sets, continuity of $V^*(x)$ is guaranteed by \cite[Theorem C.34]{MPC_book}. \qed 

\begin{lemma}\label{closed}
The set $\mathbb{X}^{\xi}_{MPC}$ is closed, and $\mathbb{X}_{MPC}^0$ is  compact.
\end{lemma}

\textit{Proof. } Since all the constraint sets in our problem formulation are closed, and the intersection of closed sets is closed by De Morgan's laws, the feasible sets $\mathbb{X}_{MPC}^0$ and $\mathbb{X}^{\xi}_{MPC}$ will be closed. Given input, state and terminal state constraints $\mathbb{U}, \mathbb{X}, \mathbb{X}_f$ are compact, the feasible set $X^0_{MPC}$ of the hard-constrained MPC problem is also compact. \qed

\begin{lemma}\label{V_lpd}
The value function $V^*$ is a locally positive definite (l.p.d.) function around $\mathbb{X}_{MPC}^0$ in $\mathbb{X}^{\xi}_{MPC}$.
\end{lemma}

\textit{Proof.} $V^*: \mathbb{X}^{\xi}_{MPC}\rightarrow \mathbb{R}$ is a continuous function, by Lemma III.2. For all $x(k) \in \mathbb{X}_{MPC}^0$, the state constraints $\mathbb{X}$ along the open-loop trajectory can be satisfied without any violation, meaning $ \xi^*_{i|k} = 0,\; \forall i \in \{0, 1, ..., N-1\}$, and $V^*(x(k)) = 0$. For all $x(k) \in \mathbb{X}^{\xi}_{MPC}\backslash \mathbb{X}_{MPC}^0$, the state constraints $\mathbb{X}$ can only be satisfied with violation, meaning $ \exists \, i \in \{0, 1, ..., N-1\}, s.t. \; \xi^*_{i|k} > 0$, and $V^*(x(k)) > 0$. By Definition \ref{lpd}, $V^*$ is a l.p.d. function around $\mathbb{X}^0_{MPC}$ in $\mathbb{X}^{\xi}_{MPC}$. \qed

\begin{lemma}\label{converge_O}
Problem \eqref{soft_MPC} defines a controller \( u_0^*(x) \). If \( \kappa(x) = u_0^*(x) \) in \eqref{closed_loop_system}, then the closed-loop slack variables $\xi_0^*(x(k)) = \xi^*_{0|k}$ converge to zero, as $k \rightarrow \infty$.
\end{lemma}
\textit{Proof.} Denote the optimal input, state and constraint violation sequence\footnote{In case of multiple minimizers, a solution $\textbf{u}^*_k$ can be singled out by the application of consistent tie-break rules, e.g., minimizing the $l_2$-norm of $\textbf{u}^*_k$.} for problem \eqref{soft_MPC} at current time $k$ by
\begin{equation}\label{optimal_sequence}
\begin{split}
    \textbf{u}^*_k & = \{u^*_{0|k}, u^*_{1|k}, \cdots, u^*_{N-1|k} \}, \\
    \mathbf{x}^*_k & = \{x^*_{0|k}, x^*_{1|k}, \cdots, x^*_{N|k} \},\\
    \Xi^*_k & = \{\xi^*_{0|k}, \xi^*_{1|k}, \cdots, \xi^*_{N-1|k} \}
\end{split}
\end{equation}
where $u^*_{i|k} \in \mathbb{U}, \; \xi^*_{i|k} \geq 0, x^*_{i|k} \in \mathbb{X}(\xi^*_{i|k}), \; \forall i \in \{0, 1, ..., N-1\}$ and $x^*_{N|k} \in \mathbb{X}_f$. The corresponding value function is
\begin{equation*}
V^*(x(k)) = \sum_{i=0}^{N-1} \xi^*_{i|k}.
\end{equation*}
A feasible solution at the $k+1^{th}$ step is 
\begin{equation}\label{feasible_sequence}
\begin{split}
    \hat{\textbf{u}}_{k+1} & = \{u^*_{1|k}, u^*_{2|k}, \cdots, u^*_{N-1|k}, \kappa_f(x^*_{N|k}) \}, \\
    \hat{\mathbf{x}}_{k+1} & = \{x^*_{1|k}, x^*_{2|k}, \cdots,x^*_{N-1|k}, x^*_{N|k}, f(x^*_{N|k}, 
    \kappa_f(x^*_{N|k})) \},\\
    \hat{\Xi}_{k+1} & = \{\xi^*_{1|k}, \xi^*_{2|k}, \cdots, \xi^*_{N-1|k}, 0\}.
\end{split}
\end{equation}
We need to verify if the proposed solution is feasible by checking if it satisfies the hard constraints $\mathbb{U} $ and $ \mathbb{X}_f$. Since there exists $\kappa_f(x^*_{N|k}) \in \mathbb{U}$, given $x^*_{N|k} \in \mathbb{X}_f$, the proposed control sequence $\hat{\textbf{u}}_{k+1}$ satisfies all the input constraints. Since $\mathbb{X}_f$ is positively invariant under $\kappa_f$, we have $f(x^*_{N|k}, \kappa_f(x^*_{N|k})) \in \mathbb{X}_f$, so the proposed state sequence $\hat{\mathbf{x}}_{k+1}$ satisfies all the state constraints. Since $x^*_{N|k} \in \mathbb{X}_f \subset \mathbb{X}(\hat{\xi}_{N-1|k+1} = 0)$, we conclude that the proposed solution \eqref{feasible_sequence} is feasible for problem \eqref{soft_MPC}.

The value function at the $k+1^{th}$ step is bounded by  
\begin{equation}\label{decrease_property}
    V^*(x(k+1)) \leq \sum_{i=1}^{N-1} \xi^*_{i|k} = V^*(x(k)) - \xi^*_{0|k}.
\end{equation}

The inequality \eqref{decrease_property} is a specific instance of the inequality presented in Lemma \ref{convergence}, where $l(k) = V^*(x(k))$, $\zeta(k) = 0$ and $w(k) = \xi^*_{0|k}$. By Lemma \ref{convergence}, it follows that the sequence $(V^*(x(k)))_{k\in \mathbb{N}}$ converges and the sequence $(\xi^*_{0|k})_{k\in \mathbb{N}}$ is summable. Since $(\xi^*_{0|k})_{k\in \mathbb{N}}$ is a non-negative sequence, to be summable, it means $\xi^*_{0|k} \rightarrow 0$ as $k \rightarrow \infty$\footnote{For all $x(k) \in \mathbb{X}$, $\xi^*_{0|k} = 0$, thus, $\xi^*_{0|k}$ is not a l.p.d. function around $\mathbb{X}_{MPC}^0$ in $\mathbb{X}^{\xi}_{MPC}$. Therefore, the Lyapunov theorem for a set can not be used to prove asymptotic stability of safe invariant sets within $\mathbb{X}$.}. \qed

Consider a trajectory starting from an arbitrary point $x(0) \in \mathbb{X}^{\xi}_{MPC}$. The trajectory must stay in $\mathbb{X}^{\xi}_{MPC}$ all the time since $\mathbb{X}^{\xi}_{MPC}$ is positively invariant by Lemma \ref{lemma_positive_invariant}. As a result, $\xi^*_{0|k} \rightarrow 0$ implies that all the system trajectories starting from within $\mathbb{X}^{\xi}_{MPC}$ converge to the set $\mathbf{R} := \{x \in \mathbb{X}^{\xi}_{MPC}:\xi^*_{0}(x) = 0\} = \mathbb{X}^{\xi}_{MPC}\cap \mathbb{X}$. 

\begin{lemma}\label{bounded_trajectory}
Problem \eqref{soft_MPC} defines a controller \( u_0^*(x) \). If \( \kappa(x) = u_0^*(x) \) in system \eqref{closed_loop_system}, then all closed-loop trajectories, starting from $\mathbb{X}^{\xi}_{MPC}$, stay within some compact subsets of $\mathbb{X}^\xi_{MPC}$ that are defined by their initial conditions $x(0)$.
\end{lemma}

\textit{Proof.} For every feasible initial condition $x(0) \in \mathbb{X}^{\xi}_{MPC}$, there exists optimal open-loop slack sequence $\{\xi^*_{0|0}, \xi^*_{1|0}, \cdots, \xi^*_{N-1|0} \}$ and corresponding value function $V^*(x(0)) = \sum_{i=0}^{N-1} \xi^*_{i|0}$. Since the value function is non-increasing, $V^*(x(0))$ provides an upper bound for the value function of subsequent trajectory, which can be considered as an additional constraint that is imposed implicitly on the slack variables at subsequent time steps. For instance, at $k^{th}$ time step, apart from the constraints shown in \eqref{soft_MPC}, there is an implicitly encoded constraint that 
\begin{equation}\label{implicit_constraint}
\sum_{i = 0}^{N-1} \xi_{i|k} \leq \sum_{i = 0}^{N-1} \xi^*_{i|0}
\end{equation}
Since input space and state space are compact with continuous dynamics, \eqref{implicit_constraint} plus the constraints in \eqref{soft_MPC} define an invariant compact set. 
Therefore, $V^*(x(0))$ defines a compact subset of $\mathbb{X}^\xi_{MPC}$ and the closed-loop trajectory of $x(0)$ remains in the compact subset. Therefore, the closed-loop trajectory $\cup_{k=0}^\infty \{x(k)\}$ for any $x(0) \in \mathbb{X}^\xi_{MPC}$ is bounded. \qed

\begin{theorem}
Problem \eqref{soft_MPC} defines a controller \( u_0^*(x) \). If \( \kappa(x) = u_0^*(x) \) in system \eqref{closed_loop_system}, then for all $x(0) \in \mathbb{X}^{\xi}_{MPC}$, its closed-loop trajectory will converge to the maximal invariant set within $\mathbb{X}^{\xi}_{MPC} \cap \mathbb{X}$. 
\end{theorem}

\textit{Proof.} We have shown that any closed-loop trajectories in $\mathbb{X}^{\xi}_{MPC}$ eventually converges to $\mathbf{R}$ and the trajectories are bounded according to Lemma \ref{bounded_trajectory}. Any \emph{bounded} trajectory of an autonomous system converges to an invariant set, which is the positive limit set of the trajectory \cite[Section 3.4.3]{Nonlinear_Control_Book}. Thus, any trajectory will converge to a subset of the union of all invariant sets within $\mathbf{R}$. The union of all invariant sets within $\mathbf{R}$ is denoted by $\mathbf{M}$, which is the maximal invariant set within $\mathbf{R}$, i.e., $\mathbb{X}^{\xi}_{MPC}\cap \mathbb{X}$.\qed


\begin{lemma} $\mathbb{X}_{MPC}^0$ is a subset of $\mathbf{M}$. \end{lemma}

\textit{Proof.} $\mathbb{X}_{MPC}^0$ is defined as the feasible set for the hard-constrained MPC problem. Specifically, consider a trajectory, if $x(0) \in \mathbb{X}_{MPC}^0,$ it also has the property that $x(0) \in \mathbb{X}^{\xi}_{MPC}\cap \mathbb{X}$. Since $\mathbb{X}_{MPC}^0$  has proved to be positively invariant by Lemma \ref{lemma_positive_invariant}, $x(k) \in \mathbb{X}_{MPC}^0$ for all $k \geq 0$, implying that $x(k) \in \mathbb{X}^{\xi}_{MPC}\cap \mathbb{X}$, for all $k\geq 0$. Thus, $\mathbb{X}_{MPC}^0$ is an invariant set within $\mathbb{X}^{\xi}_{MPC}\cap \mathbb{X}$, as a result, $\mathbb{X}_{MPC}^0\subseteq \mathbf{M}$. \qed

\begin{theorem}
Problem \eqref{soft_MPC} defines a controller \( u_0^*(x) \). If \( \kappa(x) = u_0^*(x) \) in system \eqref{closed_loop_system}, then all its closed-loop trajectories with $x(0) \in \mathbf{M}$ satisfy state constraints throughout.
\end{theorem}

\textit{Proof.} We consider two cases. The first one is for $x(0) \in \mathbb{X}_{MPC}^0$. The second is when $x(0) \in \mathbf{M} \backslash \mathbb{X}_{MPC}^0$. 
Consider the first cases, the hard-constrained MPC problem guarantees the state constraint satisfaction recursively, thus, for all $x(0) \in \mathbb{X}_{MPC}^0$, $x(k) \in \mathbb{X}$ for all $k \geq 0$. Now consider the second case that $x(0) \in \mathbf{M} \backslash \mathbb{X}_{MPC}^0$. Since $\textbf{M}$ is the union of the invariant set within $\mathbb{X}^{\xi}_{MPC}\cap \mathbb{X}$, for all $x(0) \in \mathbf{M}$, its subsequent state $x(k) \in \mathbb{X}^{\xi}_{MPC}\cap \mathbb{X}$ for all $k \geq 0$, it implies $x(k)$ always satisfy the state constraints. Therefore, all the closed-loop  trajectories start within set $\mathbf{M}$ satisfy the state constraints and $\mathbf{M}$ is a safe invariant set. \qed

In summary, based on Theorems 1 and 2, we have shown that all trajectories that start from $\mathbb{X}^\xi_{MPC}$, with the safe MPC controller $\kappa(x) = u_{0}^*(x)$, will eventually converge to the maximal safe invariant set $\mathbf{M}$ within $\mathbb{X} \cap \mathbb{X}^{\xi}_{MPC}$, meaning that safety will be recovered.

\subsection{Stability of the MPC formulation}
The previous section proves the convergence property of the closed-loop system \eqref{closed_loop_system} in set $\mathbb{X}^{\xi}_{MPC}$, with the safe MPC controller $\kappa(x) = u^*_{0}(x)$.
This section will prove the set $\mathbb{X}_{MPC}^0$ is a stable set for the closed-loop system in $\mathbb{X}^{\xi}_{MPC}$.

\begin{theorem}
Consider problem \eqref{soft_MPC}. The set $\mathbb{X}_{MPC}^0$ is a stable set for the closed-loop system in $\mathbb{X}^{\xi}_{MPC}$.
\end{theorem}

\textit{Proof.} This proof adopts the proof in \cite[Section 3.1]{Nonlinear_System_Book}, which establishes the standard Lyapunov stability theorem around a point. The approach has been extended here to demonstrate stability for a set.

For any $\epsilon > 0$, choose $r \in \left(0, \epsilon \right]$ such that 
\begin{equation*}
    B_r = \{x \in \mathbb{R}^n|\; |x|_{\mathbb{X}_{MPC}^0} \leq r\} \subset \mathbb{X}^{\xi}_{MPC}
\end{equation*}
Let $\alpha = \underset{|x|_{\mathbb{X}_{MPC}^0} = r}{\min} V^*(x).$ By Lemma \ref{V_lpd}, we have $\alpha > 0.$ Take $\beta \in \left( 0, \alpha \right)$ and let 
\begin{equation*}
    \Omega_\beta = \{x \in B_r | V^*(x) \leq \beta\}
\end{equation*}
Then, $\Omega_\beta$ is in the interior of $B_r$, since $V^*(\cdot)$ is a l.p.d. function. The set $\Omega_\beta$ is positively invariant that any trajectory of the system starting in $\Omega_\beta$ at $t=0$ stays in $\Omega_\beta$ for all $t \geq 0$, which follows from \eqref{decrease_property} that
\begin{equation*}
\begin{split}
    & \Delta V^*(x(k)) = V^*(x(k+1))-V^*(x(k)) \leq 0 \\
    \Rightarrow \; & V^*(x(k)) \leq V^*(x(0)) \leq \beta, \; \forall k \geq 0.
\end{split}
\end{equation*}
As $V^*(x)$ is continuous and $V^*(x) = 0, \forall x \in \mathbb{X}_{MPC}^0$, there exists $\delta > 0$ such that 
\begin{equation*}
    V^*(x) < \beta: \forall x, \text{ s.t. } |x|_{\mathbb{X}_{MPC}^0} \leq \delta
\end{equation*}
Then, 
\begin{equation*}
    B_\delta \subset \Omega_\beta \subset B_r
\end{equation*}
For any $x(0) \in B_\delta$, we have $x(0) \in \Omega_\beta$. Since $\Omega_\beta$ is positively invariant, then $x(k) \in \Omega_\beta$ for all $k \geq 0$, which further implies $x(k) \in B_r$.
Therefore, for any $\epsilon > 0$ and $0 < r \leq \epsilon$, there exists $\delta > 0$ so that if $|x(0)|_{\mathbb{X}_{MPC}^0} < \delta$, then $|x(k)|_{\mathbb{X}_{MPC}^0} < r \leq \epsilon$ for all $k \geq 0$, which shows the set $\mathbb{X}_{MPC}^0$ is stable. \qed

We have separately proved convergence to the safe invariant set and the set stability property using the invariant set theorem and Lyapunov theorem. Using the invariant set theorem allows us to prove convergence with a non-increasing value function, while in \cite{Predictive_CBF}, they proved the set asymptotic stability with the Lyapunov theorem which requires a strictly decreasing value function.

\subsection{Predictive control barrier function (PCBF)}

\begin{theorem}
Consider problem \eqref{soft_MPC}. The value function $V^*$ is a Discrete-Time CBF for the safe set $\mathbb{X}_{MPC}^0$ with domain $\mathbb{X}^{\xi}_{MPC}$.
\end{theorem}

\textit{Proof.} Since the value function $V^*$ has proven to be a l.p.d. fucntion around $\mathbb{X}_{MPC}^0$ in $\mathbb{X}^{\xi}_{MPC}$ by Lemma \ref{V_lpd}, we have:
\begin{equation*}
\begin{split}
V^*(x(k)) & \left\{
\begin{array}{ll}
      =0, & \text{if } x(k) \in \mathbb{X}_{MPC}^0, \\
      >0, & \text{if } x(k) \in \mathbb{X}^{\xi}_{MPC}\backslash \mathbb{X}_{MPC}^0.
\end{array}
\right. \\
\end{split}
\end{equation*}
In our problem, $\mathbb{X}_{MPC}^0$ is the safe set, and the barrier function we aim to verify is $V^*$ with domain $\mathbb{X}^{\xi}_{MPC}$. Recalling that by \eqref{decrease_property} that $V^*(x(k+1)) -  V^*(x(k)) \leq - \xi^*_{0|k}$. For all $x \in \mathbb{X}$, we have $\xi^*_{0}(x)=0$, meanwhile for all $x \in \mathbb{X}^{\xi}_{MPC}\backslash \mathbb{X}$, it holds that $\xi^*_{0}(x)>0$. Thus, we conclude that $\inf_{u\in\mathbf{U}} [V^*(f(x,u)) - V^*(x)] \leq 0$ for all $x \in \mathbb{X}^{\xi}_{MPC}$, condition \eqref{self_defined_h_condition} is satisfied. Therefore, $V^*$ is indeed a Discrete-Time CBF for the safe set $\mathbb{X}_{MPC}^0$ with domain $\mathbb{X}^{\xi}_{MPC}$.\qed

\textbf{Remarks.} The value function of the safe MPC is a \emph{Zeroing Control Barrier Function (ZCBF)} \cite{ZCBF}, which is obtained by incorporating state constraints into the objective function using penalty functions. We conjecture that if barrier functions are used to incorporate state constraints into the objective function, the resulting value function $V^*$ will be a \emph{Reciprocal Control Barrier function (RCBF)}, but this approach does not extend the domain of the MPC, i.e., the domain of the resulting RCBF will be $\mathbb{X}_{MPC}^0$, instead of $\mathbb{X}^{\xi}_{MPC}$. This outcome aligns with the result of \cite{Robustness_CBF} that ZCBF exhibits a robustness property, whereas RCBF does not. We will investigate the relationship with RCBF in future work.


\section{Numerical Examples}
\label{sec: exmaples}
This section presents two examples\footnote{While we use 2D examples for illustrative purpose, our safe MPC formulation naturally extends to higher-dimensional problems, like standard MPC approaches.} illustrating the efficacy of our results. The first demonstrates the convergence property of the proposed MPC formulation using an unstable linear system, showing that the system can recover from initially unsafe states to safe states. The second example, which is nonlinear, compares the safe invariant set of PCBFs with competing methods. In both examples, we use the MPT3 toolbox \cite{MPT3} to compute maximal invariant sets. In the linear example, the MPC problem is solved using CVXPY \cite{CVXPY} with the Clarabel solver \cite{Clarabel_2024}. In the nonlinear example, we use the Casadi framework \cite{Casadi} with the IPOPT solver \cite{IPOPT} to solve the nonlinear MPC problem.

\subsection{Unstable Linear Example}
We consider the unstable linear system
\begin{equation}
    x(k+1) = \begin{bmatrix}
        1.5 & 1 \\
        0 & 1
    \end{bmatrix} x(k) + \begin{bmatrix}
        0.5 \\
        0.5
    \end{bmatrix}u(k)
\end{equation}
which is subject to state constraints $||x||_\infty \leq 1$ and input constraints $||u||_\infty \leq 1.5$. A local stabilizing linear quadratic regulator (LQR) controller is defined as $u(k) = K_px(k)$ with $K_p = [-1.3735, -1.6166]$. It is a stabilizing terminal control law on the terminal invariant set $\mathbb{X}_f = \{x \in \mathbb{R}^n:x^TPx \leq \alpha\}$ with $P = \begin{bmatrix}
    3.3729 & 0.3776\\
    0.3776 & 1.1956
\end{bmatrix}, \alpha = 0.6$. 
We use a planning horizon $N = 10$.

\begin{figure}[ht]
    \centering   \includegraphics[width=0.48\textwidth]{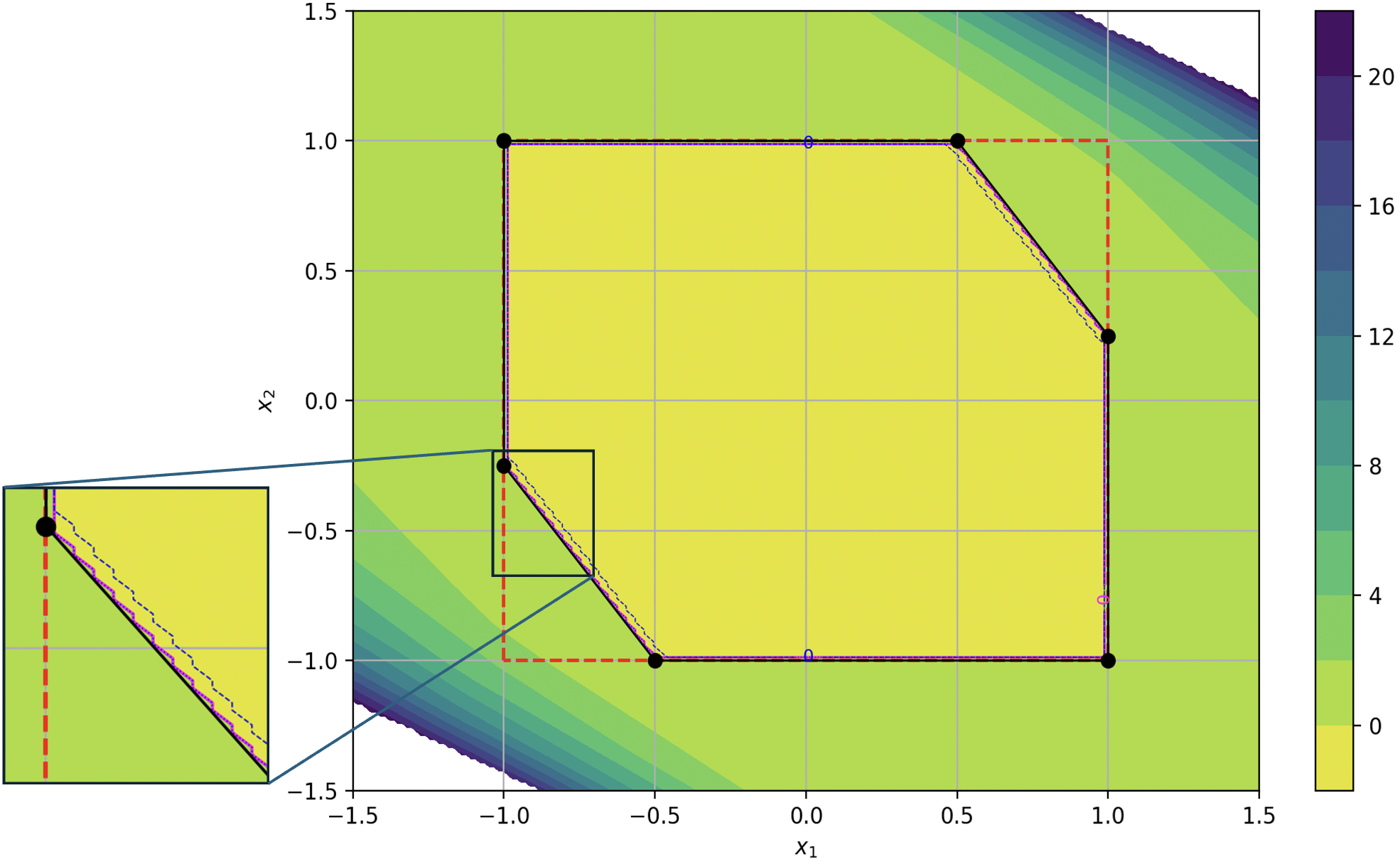}
        \caption{This is a filled contour plot of the value function $V^*$. The state space $\mathbb{X}$ is denoted with the red dashed line. The maximal safe invariant set is denoted with the black solid line. The zero-level set of $V^*$, corresponding to the predictive safe set $\mathbb{X}_{MPC}^0$, is marked with the purple solid line. The predictive safe set obtained with constraints tightening $\Delta_i = i \cdot 0.05, i \cdot 0.005$ is plotted with dashed and dotted lines, respectively. The dotted line for $\Delta_i = i \cdot 0.005$ overlaps with the zero-contour of $V^*$. }
        \vspace{-5pt}
    \label{fig:modified level set}
\end{figure}

Figure~\ref{fig:modified level set} shows the relative size of the state constraint set, the maximal safe invariant set, and the safe invariant set encoded by the PCBF $V^*$, i.e., $\mathbb{X}_{MPC}^0 = \{x\in \mathbb{R}^n: V^*(x) = 0 \}$ and other contour lines of $V^*$. The finite-value contours indicate the regions that eventually converge to the maximal safe invariant set. We also include the predictive safe set obtained by applying the MPC formulation from \cite{Predictive_CBF}, noting that this approach depends on the choice of a tightening parameters $\Delta_i$. To observe the tightening effect, we use $\Delta_i = i \cdot 0.05, i \cdot 0.005$. For the case $\Delta_i = i \cdot 0.005$, the predictive safe set overlaps with the one without tightening. As shown in Figure~\ref{fig:modified level set}, our formulation yields a predictive safe set that better approximates the maximal safe invariant set, which is expected since our method does not involve constraint tightening while still ensuring safety recovery. 
As the tightening parameters $\Delta_i$ decrease, the predictive safe set obtained converges to the set produced by our formulation, $\mathbb{X}_{MPC}^0$.

\begin{figure}[ht]
    \centering   \includegraphics[width=0.45\textwidth]{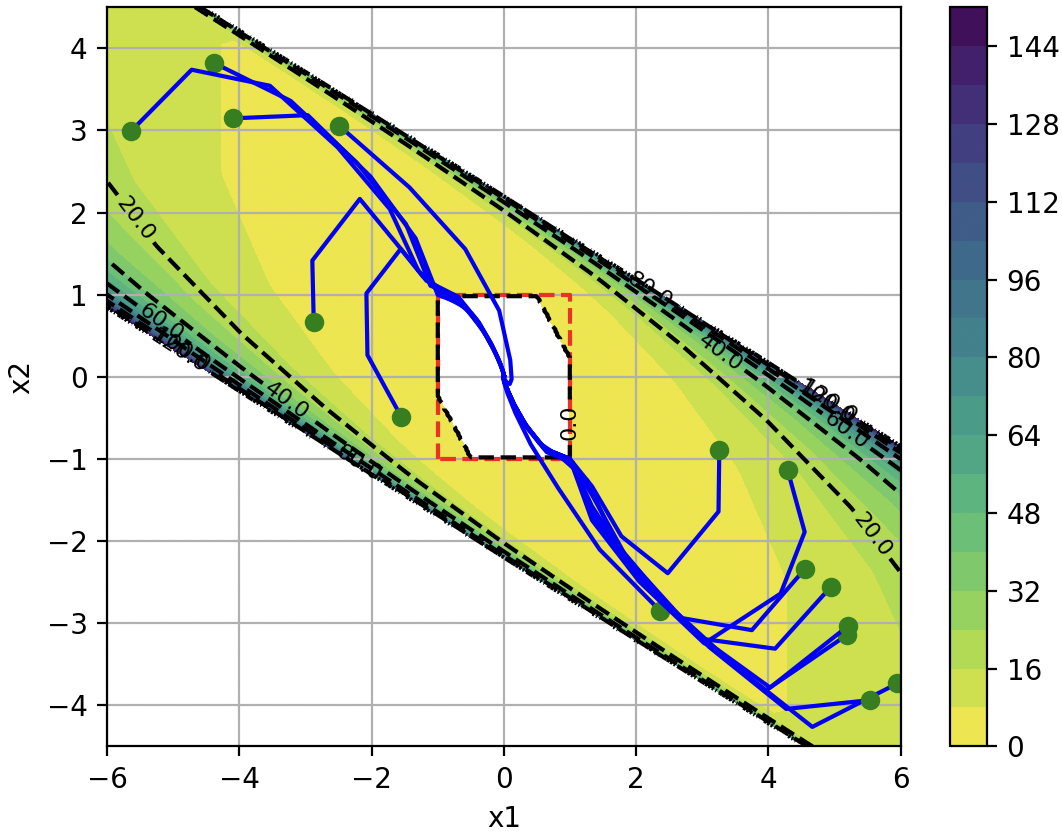}
    \vspace{-5pt}
    \caption{The safe set $\mathbb{X}$ is denoted by the red dashed line. The safe invariant set $\mathbb{X}_{MPC}^0$ is marked with the black dashed line. Closed-loop trajectories are the blue lines, they have randomly sampled initial conditions and all converge to the safe invariant set within $\mathbb{X}$, given input constraint $\mathbb{U}$.}
    \vspace{-15pt}
    \label{fig:PSF trajectory modified linear}
\end{figure}

The safe MPC formulation can serve as an auxiliary problem in conjunction with Predictive Safety Filters (PSFs), following \cite[Algorithm 1]{Predictive_CBF}, and we choose the desired control input $u_p(x) = K_px$.
This algorithm ensures safety by making minimal adjustments to the desired control input. 
By incorporating the auxiliary problem, we introduce a recovery mechanism that helps navigate out of unsafe situations as demonstrated in Fig. \ref{fig:PSF trajectory modified linear}. In Fig. \ref{fig:PSF trajectory modified linear}, we plot the safe set $\mathbb{X}$, the safe invariant set $\mathbb{X}_{MPC}^0$ and closed-loop trajectories, whose initial states are sampled randomly from unsafe regions so that $x(0) \notin \mathbb{X}$ and $u_p(x(0)) \notin \mathbb{U}$. As guaranteed by Theorems 1 and 2, all trajectories converge to the maximal safe invariant set of the system.

\subsection{Nonlinear Example}
\begin{figure}[ht]
    \centering   \includegraphics[width=0.45\textwidth]{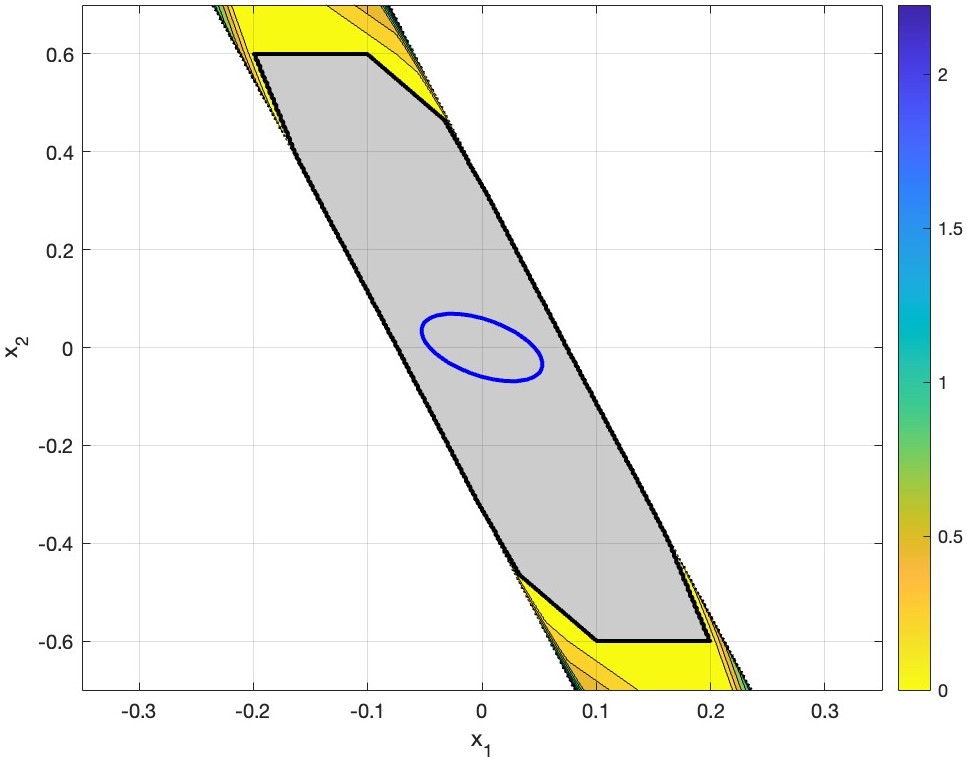}
    \vspace{-5pt}
    \caption{Compare the safe invariant sets derived from various methods. The predictive safe set of the PCBF, the maximal invariant set, and the safe invariant set of the hand-crafted CBF are represented by the black line, grey region, and blue line, respectively. It also shows the contours of PCBF $V^*$.
}
    \vspace{-10pt}
    \label{fig:nonlinear contour comparison}
\end{figure}

We next consier a nonlinear example to compare the size of invariant sets associated with CBFs generated by different methods.
Consider the nonlinear system
\begin{equation}\label{exmp:nonlinear}
    \frac{d}{dt} \begin{bmatrix}
        x_1 \\
        x_2
    \end{bmatrix} = \begin{bmatrix}
        x_2 \\
        10\sin(2x_1)
    \end{bmatrix} + \begin{bmatrix}
        0 \\
        \frac{1}{2}
    \end{bmatrix}u
\end{equation}
The system is subject to state constraints $|x_1| \leq 0.3, \; |x_2| \leq 0.6$, along with an input constraint $|u| \leq 3$. The system is discretized using the Euler forward method with a step size of $\Delta T = 0.5$, meaning $\mathbf{x}(k+1) = \mathbf{x}(k) + 0.5 \cdot \frac{d\mathbf{x}}{dt}(k) $. To compute the PCBF, following the safe MPC formulation in \eqref{soft_MPC}, the terminal set $\mathbb{X}_f$ is chosen as $(x_1, x_2) = (0, 0)$ with a prediction horizon of $N = 10$. The value function $V^*$ of the MPC is a PCBF, and its zero-level set defines the corresponding safe invariant set, represented by the black line in Fig. \ref{fig:nonlinear contour comparison}. To approximate the maximal positively invariant set, we use the MPT3 toolbox \cite{MPT3} by modelling the nonlinear dynamics as piecewise affine systems, which is depicted as the grey-filled region. 
A candidate CBF $h(\mathbf{x}) = 1- \frac{x_1^2}{a^2} - \frac{x_2^2}{b^2} - \frac{x_1 x_2}{ab}$  with parameter $a, b > 0$, is proposed in \cite{Alan2022ControlBF} for this system, based on an alternative CBF definition \cite{Alternative_CBF_definition}. 
This choice of CBF does not account for input constraints, so the parameters $a $ and $ b$ must be chosen such that input constraints are satisfied throughout the invariant set. For this example we select, $a = 0.046$ and $b = 0.06$ so that $h(\mathbf{x})$ acts as a control barrier function satisfying both state and input constraints, and it is illustrated with a blue line in Fig \ref{fig:nonlinear contour comparison}. 

The PCBF is effective in expanding the predefined invariant set. In this example, the predefined invariant set is limited to the origin, yet the resulting PCBF produces a predictive safe set that nearly matches the maximal invariant set. In contrast, constructing a non-conservative CBF that accounts for input constraints is challenging. In addition to achieving a non-conservative CBF, a key benefit of using our safe MPC approach its ability to identify a region from which recovery to the safe invariant set is possible even from unsafe states, while also providing the recovery controller. Within this recoverable region, the value function remains finite and positive. The level sets of the value function are illustrated in Figure \ref{fig:nonlinear contour comparison}. Due to input constraints, the safe MPC problem may not be globally feasible. 

\section{Conclusion}
In this paper we propose a safe MPC formulation to ensure invariance. It connects MPC techniques with CBFs, by having its value function being a CBF, termed as PCBF. The safe MPC formulation effectively expands a predefined conservative safe invariant set to give a good approximation of the maximal safe invariant set. All trajectories that start within the domain of the PCBF, including points in unsafe regions, eventually converge to the maximal safe invariant set and therefore recover safety. We prove the convergence property with the invariant set theorem that only requires the non-increasing value function. This controller aims to guarantee safety rather than stability, however, stability can also be addressed by incorporating CLF constraints \cite{ZCBF}. Future work aims at establishing a quantitative guarantee regarding the relative size of the safe invariant set of the PCBF compared to the system's maximal safe invariant set.




\bibliographystyle{unsrt}
\bibliography{mybib.bib}


\end{document}